\documentclass [12pt]  {article} 
\usepackage{amssymb}
\def \CC{{\mathbb{C}}}

\def \FFF{{\cal{F}}}

\begin{document}

\begin{center}
{\Large {\bf Variation on a Theorem by Mues and Steinmetz}}\\
\bigskip
{\tiny {\rm BY}}\\
\bigskip
{\sc Andreas Schweizer}\\
\bigskip
{\small {\rm Department of Mathematics,\\
Korea Advanced Institute of Science and Technology (KAIST),\\ 
Daejeon 305-701\\
South Korea\\
e-mail: schweizer@kaist.ac.kr}}
\end{center}
\begin{abstract}
\noindent
Let $f$ be a meromorphic function. We suggest a generalization 
of $f$ and its derivative $f'$ sharing a nonzero value $a$ IM 
that does not impose any a priori restrictions on the ramification 
of $f$. Then we discuss some results around the question whether
the famous theorem on entire functions $f$ that share two 
values IM with $f'$ still holds for this weaker notion.
\end{abstract}

\subsection*{1. Introduction}

Two meromorphic functions $f$ and $g$ on a complex domain $D$ are 
said to share the value $a\in\CC$ if for all $z\in D$ we have
$f(z)=a\Leftrightarrow g(z)=a$.
Usually this is more loosely written as 
$f=a\Leftrightarrow g=a$.
Often one also says more precisely that $f$ and $g$ are sharing 
the value $a$ IM (ignoring multiplicities) to emphasize that one
does not require that at such points the function $f$ takes the
value $a$ with the same multiplicity as $g$.
\par
Here we are interested in the case where the domain $D$ is the
whole complex plane $\CC$, $f$ is an entire function, and $g$
is its derivative $f'$. The most famous theorem in this setting 
is the following result from 1979.
\\ \\
{\bf Theorem A.} [MuSt, Satz 1] \it
Let $f$ be a nonconstant entire function and let $a$ and $b$ be complex 
numbers with $a\neq b$. If $f$ and $f'$ share the values $a$ and $b$ IM, 
then $f\equiv f'$.
\rm
\\ \\
See also [YaYi, Theorem 8.3] for a proof in English. A different proof 
for the case $ab\neq 0$ was also given in [Gu, Theorem 2].
\\ \\
Theorem A has been generalized in many ways. We will use 
a relatively new result which gives up one of the directions
of the implication $f=b\Leftrightarrow f'=b$. Actually, it 
also provides a third different proof for the case 
$ab\neq 0$ of Theorem A.
\\ \\
{\bf Theorem B.} [L\"uXuYi, Theorem 1.1] \it
Let $a$ and $b$ be two nonzero distinct complex numbers, and let 
$f(z)$ be a nonconstant entire function. If 
$f(z)=a\Leftrightarrow f'(z)=a$ and 
$f'(z)=b\Rightarrow f(z)=b$, then $f\equiv f'$.
\rm
\\ \\
The functions $f=\frac{b}{4}z^2$ and $f=Ce^{\frac{a}{C}z}+a-C$
from [LiYi] show that the conditions $a\neq 0$ and $b\neq 0$
cannot be removed in Theorem B. Here $C$ is any nonzero constant.
\par
There is also a result that gives up the other direction of the 
implication. It provides an alternative proof for the case
$ab=0$ of Theorem A.
\\ \\
{\bf Theorem C.} [L\"uXuYi, special case of Corollary 1.1] \it
Let $b$ be a nonzero number, and let $f(z)$ be a nonconstant entire 
function. If $f(z)=0\Rightarrow f'(z)=0$, and $f$ and $f'$ share
$b$ IM, then $f\equiv f'$.
\rm
\\ \\
Actually, [L\"uXuYi, Corollary 1.1] states more generally 
that $f=0\Rightarrow f'=0$ and $f=b\Rightarrow f'=b$ implies 
$f\equiv f'$ or $f=b(\frac{A}{2}e^{\frac{z}{4}}+1)^2$ with 
a nonzero constant $A$. The additional condition 
$f'=b\Rightarrow f=b$ excludes the second possibility.
\par
A more general version of Theorem C allowing $a\neq 0$ 
(compare [L\"uXuYi, Theorem 1.2]) would also admit 
solutions with $f\not\equiv f'$, 
namely $f=Ce^{\frac{b}{b-a}z}+a$ for $b\neq 0$, and 
$f=a(\frac{A}{2}e^{\frac{z}{4}}+1)^2$ for $b=0$.
\\ \\
The starting point of our paper is the following observation.
If a meromorphic function $f$ shares the nonzero value $b$ with 
its derivative $f'$, then, in contrast to $f$ sharing $b$ with 
an arbitrary meromorphic function $g$, this imposes the additional 
condition that all $b$-points of $f$ are simple. (This is actually 
used in the proofs of all theorems above.) So one might wonder 
about value sharing between $f$ and $f'$ that does not impose 
any conditions on the ramification of $f$. 
\par
In one direction, the best that we can demand is that
$f=b$ implies $f'=b$ unless for ramification reasons we
have $f'=0$. In other words, every simple $b$-point of
$f$ is a (not necessarily simple) $b$-point of $f'$.
On the other hand we can still insist on
$f'=b\Rightarrow f=b$, as this does not say anything
about possible multiple $b$-points of $f$.
\\ \\
{\bf Definition.}
Let $f$  be a meromorphic function on some domain $D$ and 
$b$ a nonzero complex number. We say that
$f$ and $f'$ {\bf share the value $b$ allowing ramification} 
if for all $z\in D$ we have 
$$f(z)=b\Rightarrow f'(z)\in\{b,0\}\ \ \ \hbox{\rm and}\ \ \ \ 
f'(z)=b\Rightarrow f(z)=b.$$
In other words, the $b$-points of $f'$ are exactly the 
simple $b$-points of $f$.
\par
In some other notation this would be written as
$$\overline{E}(b,f')=\overline{E}_{1)}(b,f).$$

Obviously it is a generalization of $f$ and $f'$ sharing the value 
$b$ IM. But in general it is neither weaker nor stronger than the 
notion of $f$ and $f'$ sharing their simple $b$-points, which we 
have discussed a bit in [Sch]. 
\par
A priori, this sharing does not impose any restrictions on the
ramification of $f$ nor of $f'$ (beyond the general condition
that outside the poles $f$ is ramified exactly at the zeroes 
of $f'$, which already holds without the sharing).
\par
So it might come as a surprise that, at least for entire functions,
sharing $b$ allowing ramification actually does impose {\it some} 
restriction on the ramification.
\\ \\
{\bf Lemma 1.1.} \it
Let $f$ be a nonconstant entire function and $b$ a nonzero complex 
number. If $f$ and $f'$ share the value $b$ allowing ramification,
then $b$ cannot be a totally ramified value of $f$.
\rm
\\ \\
This will be part of Theorem 2.6.
\par
For the rest of the paper we will only consider nonconstant entire
functions $f$. If $f$ and $f'$ share the nonzero value $b$ allowing
ramification, then by Lemma 1.1 the extremal case 
$f=b\Rightarrow f'=0$
cannot occur. The other extremal case, $f=b\Rightarrow f'=b$ simply
means that $f$ and $f'$ share the value $b$ in the usual sense (IM).
\par
Several natural questions immediately come to ones mind. However,
they seem to be less trivial than one might first think.
\\ \\
{\bf Question 1.}
Give an example of a nonconstant entire function $f$ and a nonzero 
value $b$ such that $f$ and $f'$ share $b$ allowing ramification,
but they do {\it not} share $b$ in the usual sense. In other words,
give an example of sharing $b$ allowing ramification such that at
least one of the $b$-points of $f$ really is multiple.
\\ \\
{\bf Question 2.}
Can a nonconstant polynomial and its derivative share a nonzero
value allowing ramification?
\\ \\
{\bf Question 3.}
Let $a$, $b$ be two distinct nonzero complex numbers. Let $f(z)$ be 
a nonconstant entire function. Assume that $f$ and its derivative 
$f'$ share the value $a$ allowing ramification and also share the 
value $b$ allowing ramification. Does this imply $f\equiv f'$?
\\ \\
In the next section we will give some partial results towards
Question 3, and also show that under different stronger conditions
the answer is positive.
\\

\subsection*{2. The main results}

In this section we elaborate on the question how far Theorem A,
and also Theorems B and C still hold if we weaken the sharing 
IM to sharing allowing ramification.
\par
If in Theorem A we only weaken one of the two sharings, 
the answer is easy.
\\ \\
{\bf Theorem 2.1.} \it
If the nonconstant entire function $f$ and its derivative $f'$ share the 
value $a$ IM and share the value $b$ ($\neq 0, a$) allowing ramification,
then $f\equiv f'$. 
\rm
\\ \\
If we weaken both sharings in Theorem A, we don't know the answer,
and we cannot even offer an educated guess what will happen. 
As a consolation we prove the following result.
\\ \\
{\bf Theorem 2.2.} \it
Let $a_1$, $a_2$, $a_3$ be three disctinct nonzero complex numbers.
Let $f(z)$ be a nonconstant entire function. If $f$ and its derivative 
$f'$ share the value $a_j$ allowing ramification for $j=1,2,3$, then
$f\equiv f'$.
\rm 
\\ \\
{\bf Remark.} Actually, the conditions in Theorem 2.2 are unnecessarily 
strong. As the proof will show, sharing $a_1$ allowing ramification 
and $f=a_j \Rightarrow f'\in\{a_j ,0\}$ for $j=2,3$ already implies 
$f\equiv f'$.
\\ \\
Next we make a first step towards Theorems B and A with weakened sharing.
\\ \\
{\bf Proposition 2.3.} \it
Let $a$ and $b$ be two nonzero distinct complex numbers, 
and let $f(z)$ be a nonconstant entire function. If 
$f$ and $f'$ share the value $a$ allowing ramification,
and $f'(z)=b\Rightarrow f(z)=b$, then $f$ is a transcendental 
function of order at most one.
\rm
\\ \\
{\bf Corollary 2.4.} \it
Let $f(z)$ be a nonconstant entire function. If $f$ and $f'$ share 
two different nonzero values allowing ramification, then $f$ is 
a transcendental function of order at most one.
\rm
\\ \\
Once we know that the order of $f$ is at most $1$, one might
try to prove $f\equiv f'$ by the method from [L\"uXuYi]
(which has also successfully been applied in [Sch]). But 
in the current setting we couldn't get that working. 
However, we extract some more information from the 
condition that the order is at most one.
\\ \\
{\bf Lemma 2.5.} \it
Let $f$ be a transcendental entire function of order at most one,
and let $b$ be a nonzero complex number.
\begin{itemize}
\item[a)] If $S$ is a finite subset of $\CC$, and 
$f=b\Rightarrow f'\in S$, then the derivative
$f'$ takes every nonzero value infinitely often.
\item[b)] If $f'=b\Rightarrow f=b$, and $f$ takes the value
$a$ only finitely often, then $a\neq b$ and 
$$f(z)=Ce^{\frac{b}{b-a}z}+a$$
where $C$ is a nonzero constant.
\end{itemize}
\rm
With the same method one could also prove statements about
the values of $f'$ in part b), and about the form of $f$ 
in part a). We have contented ourselves with the results 
we need.
\\ \\
{\bf Theorem 2.6} \it
Let $f$ be a nonconstant entire function such that $f$ and 
$f'$ share the nonzero value $b$ allowing ramification. Then 
the derivative $f'$ has no nonzero Picard value. 
\par
In particular, $b$ cannot be a Picard value of $f'$. So $b$ 
cannot be a totally ramified value of $f$. A fortiori, $b$ 
cannot be a Picard value of $f$. Actually, if $f$ is transcendental, 
then it must take the value $b$ infinitely often.
\rm
\\ \\
Now we weaken the sharing in Theorem C. We formulate 
the result slightly more generally.
\\ \\
{\bf Proposition 2.7.} \it
Let $f$ be a nonconstant entire function such that $f$ and 
$f'$ share the nonzero value $b$ allowing ramification.
If $f$ has a totally ramified value $a$, then $a\neq b$,
and the order of $f$ is at most one.
\rm
\\ \\
In the special case where $a$ is a Picard value of $f$ in
Proposition 2.7, we can actually nail down the function,
even under a weaker sharing condition for $b$.
\\ \\
{\bf Theorem 2.8.} \it 
Let $f$ be a transcendental entire function with
$f'=b\Rightarrow f=b$ for some nonzero value $b$.
If $f$ has a generalized Picard value $a$ (that is, 
a value that is taken only finitely often by $f$), 
then $a\neq b$ and
$$f(z)=Ce^{\frac{b}{b-a}z}+a$$
where $C$ is a nonzero constant. In particular, 
if $a=0$ we have $f\equiv f'$. 
\rm
\\ \\
If $a$ in Theorem 2.8 is a true Picard value, the 
condition that $f$ is transcendental is of course 
automatic.
\\

\subsection*{3. Some normality criteria}

Most results in this paper are proved using normal family
arguments. In this section we provide the necessary tools.
\\ \\
{\bf Theorem 3.1.} [LiYi, Theorem 3] \it
Let $\FFF$ be a family of holomorphic functions in a domain 
$D$, and let $a$ and $b$ be two finite complex numbers such that
$b\neq a, 0$. If, for each $f\in\FFF$ and $z\in D$,
$f=a\Rightarrow f'=a$, and $f'=b\Rightarrow f=b$, 
then $\FFF$ is normal in $D$.
\rm
\\ \\
{\bf Corollary 3.2.} [LiYi] \it
Let $f$ be an entire function and $a$ and $b$ two complex 
numbers with $b\neq a, 0$. If $f=a\Rightarrow f'=a$, 
and $f'=b\Rightarrow f=b$, then $f$ has order at most one.
\rm
\\ \\
For the proofs of Theorems 2.3 and 2.8 we have to generalize this.
\\ \\
{\bf Theorem 3.3.} \it
Let $\FFF$ be a family of holomorphic functions in a domain 
$D$. Let $a$ and $b$ be two finite complex numbers such that
$b\neq a, 0$, and let $S$ be a finite subset of $\CC$.
If, for each $f\in\FFF$ and $z\in D$,
$f=a\Rightarrow f'\in S$, and $f'=b\Rightarrow f=b$, 
then $\FFF$ is normal in $D$.
\rm
\\ \\
{\bf Proof.} \rm 
Obviously the conditions imply that we can always assume
$b\not\in S$.
\par
With the stronger condition $f=a\Rightarrow f'=a$ this is 
Theorem 3 from [LiYi], and our proof follows their proof 
closely. So we will be a bit sketchy and mainly emphasize
the points where one has to be careful with the more general 
conditions.
\par
Setting $h(z)=f(z)-a$ we have $|h'(z)|\leq M+1$ when $h(z)=0$
where $M$ is the maximum of the absolute values of the elements
of $S$. We can assume that $D$ is the unit disk and that the
family $\{f(z)-a\ :\ f\in\FFF\}$ is not normal at $0$.
\par
To bring this to a contradiction, as in [LiYi], we use a version
of the Pang-Zalcman Lemma [PaZa, Lemma 2], namely the one with
$$g_n(\xi)=\rho_n^{-1}\{f_n(z_n+\rho_n\xi)-a\}\to g(\xi)$$
locally uniformly with respect to the spherical metric on $\CC$,
where $g(\xi)$ is a nonconstant entire function satisfying
$$g^\#(\xi)\leq g^\#(0)=M+2.$$
First we prove $g=0\Rightarrow g'\in S$. If $g(\xi_0)=0$, by
Hurwitz's Theorem there exist $\xi_n\to\xi_0$ as $n\to\infty$
such that for sufficiently large $n$ 
$$g_n(\xi_n)=\rho_n^{-1}\{f_n(z_n+\rho_n\xi_n)-a\}=0.$$
So $f_n(z_n+\rho_n\xi_n)=a$ and hence
$g_n'(\xi_n)=f_n'(z_n+\rho_n\xi_n)\in S$. Since 
$g_n'(\xi_n)\to g'(\xi_0)$ and $S$ is discrete, we get 
$g'(\xi_0)\in S$.
\par
Next we prove that $g'(\xi)\neq b$ on $\CC$. If 
$g'(\xi)\equiv b$, then $g(\xi)=b\xi +c$, which together with 
$g=0\Rightarrow g'\in S$ contradicts $b\not\in S$. Thus
$g'(\xi)\not\equiv b$, and we can argue exactly as in [LiYi] 
to get $g'(\xi)\neq b$.
\par
Since $g'(\xi)$ is of order at most one, we now have
$g'(\xi)=b+e^{b_0 +b_1 \xi}$ with constants $b_0,b_1$.
\par
If $b_1\neq 0$, as in [LiYi, p.56] we get a contradiction,
recalling that $g=0\Rightarrow g'\in S$ and 
$\{|s-b|\ :\ s\in S\}$ is bounded.
\par
If $b_1 =0$, we have $g(\xi)=(b+e^{b_0})\xi +c_0$. From 
$g=0\Rightarrow g'\in S$ we get $b+e^{b_0}\in S$, and
hence the contradiction $g^\#(0)<M+2$.
\hfill$\Box$
\\ \\
{\bf Corollary 3.4.} \it
Let $a$ and $b$ be two finite complex numbers such that
$b\neq a, 0$, and let $S$ be a finite subset of $\CC$.
If $f$ is a transcendental entire function with
$f=a\Rightarrow f'\in S$, and $f'=b\Rightarrow f=b$, then 
$f$ has order at most one.
\rm
\\ \\
{\bf Proof.} \rm 
This is a standard conclusion that has been applied in
countless papers. By Theorem 3.3 the family 
$\{f(z+\omega)\ :\ \omega\in\CC\}$ 
is normal on $\CC$, that is, $f$ is a normal entire 
function and hence of order at most $1$. See for example 
[Mi, p.198] and [Mi, p.211] for a proof of the last 
implication.
\hfill$\Box$
\\ \\
{\bf Remark.} It seems that, with practically the same proof,
Theorem 3.3 and Corollary 3.4 still hold under the much weaker 
condition that $S$ is a closed, bounded set that does not 
contain $b$.
\\ \\
We also mention the following result which can be applied 
to some of the situations we are interested in.
\\ \\
{\bf Theorem 3.5.} [Li, Theorem 2.3] \it
Let $\FFF$ be a family of analytic functions in a domain $D$
of the complex plane, $a$, $b$, $c$, $d$ complex numbers with
$a\neq b$, and $M$ a positive number. Suppose that for all 
$f\in\FFF$
\begin{itemize}
\item[(i)] $|f'(z)|\leq M$ for 
$z\in (f-a)^{-1}(0)\cup (f-b)^{-1}(0)$; and
\item[(ii)] $(f'-c)^{-1}(0)\subseteq (f-d)^{-1}(0)$.
\end{itemize}
Then $\FFF$ is a normal family in $D$.
\rm
\\ \\
We also need a result that apparently goes back at least
to Milloux (see [Ha, page 9]).
\\ \\
{\bf Theorem 3.6.} \it
Let $f$ be a transcendental entire function. If $f'$ takes
the nonzero value $b$ only finitely often, then $f$ must 
take every value infinitely often.
\rm
\\

\subsection*{4. Proofs of the main results}

The first one is easy.
\\ \\
{\bf Proof of Theorem 2.1.} \rm 
If $a=0$, then $f'=0\Rightarrow f=0$ implies that $b$ is
actually also shared IM; so we are in the situation of 
Theorem A.
\par
If $a\neq 0$, Theorem 2.1 is just a special case of the 
more general Theorem B.
\hfill$\Box$
\\ \\
The main work for Theorem 2.2 is in the following recent result.
\\ \\
{\bf Theorem 4.1.} [Sch, Theorem 2.4] \it
Let $a_1$, $a_2$, $a_3$ be three distinct complex numbers. Nonconstant 
entire functions $f$ with $f\not\equiv f'$ and 
$$f=a_j\Rightarrow f'\in\{a_j,0\}$$
for $j=1,2,3$ exist if and only if $a_j=\zeta^ja_3$ with
$\zeta$ being a third root of unity, that is, if $(X-a_1)(X-a_2)(X-a_3)$
is of the form $X^3-\delta$.
\par
Moreover, functions with this property necessarily are of the form
$$f(z)=\frac{4\delta}{27\beta^2}e^{\frac{2}{3}z}+\beta e^{-\frac{1}{3}z}$$
with a nonzero constant $\beta$.
\par
Conversely, every function of this form has the stronger property that
every simple $a_j$-point of $f$ is a simple $a_j$-point of $f'$ for
$j=1,2,3$.
\rm
\\ \\
\rm 
{\bf Remark.} Note that there is an unfortunate typo in [Sch] in 
the formulation of Theorem 2.4 as well as towards the end of its 
proof; the term $e^{\frac{2}{3}z}$ is on both occasions given incorrectly
as $e^{\frac{3}{2}z}$. 
\\ \\
{\bf Proof of Theorem 2.2.} \rm 
If $f\not\equiv f'$, then by Theorem 4.1 we must have 
$a_j=\zeta^j c$ and $f=\frac{4c^3}{27}t^2 +\frac{1}{t}$ 
where $t=\frac{1}{\beta}e^{\frac{1}{3}z}$. Correspondingly, 
$f'=\frac{8c^3}{81}t^2 -\frac{1}{3t}$.
In particular, if $t=\frac{3(2+\sqrt{6})}{4c}$, then 
$f'=c$, one of the three shared values, 
but $f=(\sqrt{6}-\frac{1}{2})c\neq c$.
\hfill$\Box$
\\ \\
For the proof of Proposition 2.3 we need the following 
auxiliary result.
\\ \\
{\bf Lemma 4.2.} \it 
Let $a$, $b$ be two distinct complex numbers. The problem
$$f'=a\Rightarrow f=a\ \ \ \hbox{\it and } \ \ \ 
f'=b\Rightarrow f=b \qquad\qquad\qquad\qquad\qquad\qquad (*)$$
has no polynomial solutions of degree bigger than $2$.
\rm
\\ \\
{\bf Proof.} \rm
If $f$ is entire and $(*)$ holds, the auxiliary function
$$h=\frac{(f-f')f''}{(f'-a)(f'-b)}$$
is entire. If moreover $f$ is a polynomial of degree
at least $3$, then for degree reasons $h$ is a nonzero
constant. After scaling everything we can assume
$f=z^n +cz^{n-1}+\cdots$ with $n\geq 3$. Calculating the 
two highest terms of the numerator and denominator of 
$h$, we obtain the contradiction that
$$\frac{n(n-1)z^{2n-2}+(n-1)(2(n-1)c-n^2)z^{2n-3}+\cdots}
{n^2z^{2n-2}+2n(n-1)cz^{2n-3}+\cdots}$$
is constant.
\hfill$\Box$
\\ \\
{\bf Remark.} \rm
For $a+b\neq 0$, the quadratic polynomial 
$f(z)=\frac{a+b}{4}z^2 +\frac{ab}{a+b}$
is a solution of $(*)$. But it seems to be a nontrivial problem 
whether in the case $a+b=0$ there are entire solutions of $(*)$ 
other than the obvious ones $f\equiv f'$ or $f'$ a constant 
different from $a$ and $b$.
\\ \\
{\bf Proof of Proposition 2.3.} \rm
Polynomials of degree $1$ and $2$ are easily excluded. So
$f$ must be transcendental by Lemma 4.2. By Corollary 3.4 
the order of $f$ is bounded by $1$.
\hfill $\Box$
\\ \\
One can also prove Corollary 2.4 directly. For that one can manage 
without Corollary 3.4, as, due to the stronger conditions, one can 
get normality and the bound on the order from Theorem 3.5.
\\ \\
{\bf Proof of Lemma 2.5.} \rm
a) Asssume that $f'$ takes the nonzero value $a$ only finitely 
often. Since by [YaYi, Theorem 1.21] the order of $f'$ is also 
at most $1$, the Hadamard Factorization Theorem 
[YaYi, Theorem 2.5] implies
$$f'(z)=a+P(z)e^{\lambda z}$$
with $\lambda\neq 0$ and a polynomial $P(z)$. Thus 
$$f(z)=B+az+Q(z)e^{\lambda z}$$
with a constant $B$ and a polynomial $Q(z)$ with 
$P=Q'+\lambda Q$. By Theorem 3.6 there are infinitely 
many $z_0$ with $f(z_0)=b$. For these we must have 
$$e^{\lambda z_0}=\frac{b-B-az_0}{Q(z_0)}.$$ 
Moreover, 
$$f'(z_0)=
a+\lambda b-\lambda B-\lambda az_0+\frac{(b-B-az_0)Q'(z_0)}{Q(z_0)}$$
lies in the finite set $S$ for every such $z_0$.
This is only possible if the rational function
$$-\lambda az+\frac{(b-B-az)Q'(z)}{Q(z)}
=\frac{-\lambda azQ(z)+(b-B-az)Q'(z)}{Q(z)}$$
is constant, which cannot be for degree reasons.
\par
b) If $f$ takes the value $a$ only finitely often, 
then by the Hadamard Factorization Theorem
$$f(z)=a+P(z)e^{\lambda z}$$ 
with $\lambda\neq 0$ and a polynomial $P(z)$.
Theorem 3.6 together with $f'=b\Rightarrow f=b$
implies $a\neq b$ and that there are infinitely 
many $z_0$ with $f'(z_0)=b$. Then an argumentation 
similar to the one in part a) shows that $P$ must 
be a constant and $\lambda=\frac{b}{b-a}$.
\hfill $\Box$
\\ \\
{\bf Proof of Theorem 2.6.} \rm
The very last claim is an immediate consequence of 
$f'=b\Rightarrow f=b$ and Theorem 3.6.
\par
Now suppose that $a$ is a nonzero Picard value of $f'$. We 
only have show that this implies that $f$ has order at most 
$1$. Then Lemma 2.5 furnishes the desired contradiction.
\par
If $a\neq b$, we can to that purpose apply Corollary 3.2
(with the roles of $a$ and $b$ interchanged), as we vacuously 
have $f'=a\Rightarrow f=a$. 
\par
If $a=b$, we have $f=b\Rightarrow f'=0$ and vacuously
$f'=b\Rightarrow f=2b$. Hence the function $h(z)=f(z)-b$
satifies the conditions of Corollary 3.2 and thus has 
order at most $1$.
\hfill $\Box$
\\ \\
{\bf Proof of Proposition 2.7.} \rm
By Theorem 2.6 we have $a\neq b$. So we can apply
Corollary 3.4 to bound the order.
\par
Alternatively, we could get normality and order at most $1$ from 
Theorem 3.5. Or we could consider $h(z)=f(\frac{b-a}{b}z)-a$. 
Then all zeroes of $h$ are multiple, and $h$ and its derivative 
share the value $b-a$ allowing ramification. So Corollary 3.2
suffices to bound the order.
\hfill $\Box$
\\ \\
{\bf Proof of Theorem 2.8.} \rm
Theorem 3.6 together with $f'=b\Rightarrow f=b$ 
immediately implies $a\neq b$. Moreover, we have of 
course $f=a\Rightarrow f'\in S$ with a finite set 
$S$. So by Corollary 3.4 the order of $f$ is at 
most $1$. Then Lemma 2.5 takes care of the rest.
\hfill $\Box$
\\ \\
As an afterthought we mention that in final instance the proofs 
of almost all results use one version or another of the famous
Zalcman Lemma.

\subsection*{\hspace*{10.5em} References}
\begin{itemize}

\item[{[Gu]}] G.~Gundersen: \rm Meromorphic functions that share
finite values with their derivative, \it J. Math. Anal. Appl. 
\bf 75 \rm (1980), 441-446\\
(Correction: \it J. Math. Anal. Appl. \bf 86 \rm (1982), 307)

\item[{[Ha]}] W.K.~Hayman: \rm Picard values of meromorphic functions 
and their derivatives, \it Ann. of Math. \bf 70 \rm (1959), 9-42

\item[{[Li]}] B.Q.~Li: \rm On the Bloch Heuristic Principle,
normality, and totally ramified values, \it Arch. Math. (Basel)
\bf 87 \rm (2006), 417-426

\item[{[LiYi]}] J.~Li and H.~Yi: \rm Normal families and uniqueness
of entire functions and their derivatives, \it Arch. Math. (Basel)
\bf 87 \rm (2006), 52-59

\item[{[L\"uXuYi]}] F.~L\"u, J.~Xu and H. Yi: \rm Uniqueness 
theorems and normal families of entire functions and their derivatives, 
\it Ann. Polon. Math. \bf 95.1 \rm (2009), 67-75

\item[{[Mi]}] D.~Minda: \rm Yosida functions, in: \it Lectures on
Complex Analysis, Xian 1987, (C.T.~Chuang, ed.) \rm World Scientific,
Singapore, 1988, pp.197-213

\item[{[MuSt]}] E.~Mues and N.~Steinmetz: \rm Meromorphe Funktionen, 
die mit ihrer Ableitung Werte teilen, \it Manuscripta Math. \bf 29 \rm 
(1979), 195-206

\item[{[PaZa]}] X.~Pang and L.~Zalcman: \rm Normal families and shared 
values, \it Bull. London Math. Soc. \bf 32 \rm (2000), 325-331

\item[{[Sch]}] A.~Schweizer: \rm Entire functions sharing simple
$a$-points with their first derivative, \it Houston J. Math.,
to appear\rm

\item[{[YaYi]}] C.-C.~Yang and H.-X.~Yi: \it Uniqueness theory
of meromorphic functions, \rm Kluwer Academic Publishers Group,
Dordrecht, 2003

\end{itemize}

\end{document}